\documentclass[12pt,reqno]{amsart}

\usepackage{cases}
\usepackage{graphicx}
\usepackage{amssymb}
\usepackage{mathrsfs}
\usepackage{ulem}
\usepackage{xcolor}
\usepackage[all,pdf]{xy}

\newtheorem{theorem}{Theorem}[section]

\newtheorem{lemma}[theorem]{Lemma}
\newtheorem{proposition}[theorem]{Proposition}

\numberwithin{equation}{section}

\newcommand{\beqm}{\begin{eqnarray*}}
\newcommand{\eqm}{\end{eqnarray*}}

\newcommand{\C}{\mathbb{C}}

\newcommand{\abs}[1]{\left| #1 \right|}
\newcommand{\norm}[1]{\left\| #1 \right\|}

\newcommand{\whmu}{\widehat{\mu}}

\begin{document}
\title[summing Hankel operators and Berger-Coburn phenomenon]{absolutely summing Hankel operators on Fock spaces and the Berger-Coburn phenomenon}

\author[Zhangjian Hu]{Zhangjian Hu}
\address{Zhangjian Hu \\Department of Mathematics\\ Huzhou University, Huzhou
313000,   China } \email{huzj@zjhu.edu.cn}

\author[Xiaofen Lv]{Xiaofen Lv}
\address{Xiaofen Lv \\Department of Mathematics\\ Huzhou University, Huzhou
313000,   China } \email{lvxf@zjhu.edu.cn}

\thanks {Zhangjian Hu is supported by the National Natural Science Foundation of China
(Grant No. 12171150); Xiaofen Lv is supported by the National Natural Science Foundation of China
(Grant No. 12171150, 12431005)  and the National Natural Science Foundation of
Zhejiang Province (Grant No. LMS26A010001).}



\begin{abstract}

In this paper, for $1 \leq p, r < \infty$ we characterize those symbols $f$ so that the induced Hankel operators $H_f$ are   $r$-summing  from  Fock spaces $F^p_\alpha$ to  $L^p_\alpha$.  The main result shows that the $r$-summing norm of $H_f$ is equivalent to the 
 $\mathrm{IDA}^{\kappa, p}$-norm of $f$, where $\kappa$ is a positive number  determined by   $p$ and $r$, and the $\mathrm{IDA}$ space is as in \cite{HV23}.
 As some application, we discuss the  Berger-Coburn phenomenon for $r$-summing Hankel operators on  Fock spaces.

\end{abstract}

\maketitle

\section{Introduction}

The study of function spaces and associated operators has long been a central theme in complex analysis and operator theory. Among these, the Fock space, also known as the Segal-Bargmann space, plays a particularly important role due to their deep connections with quantum mechanics, time-frequency analysis, and complex geometry.

For $0<p<\infty$, $\alpha > 0$, the space $L^{p}_\alpha$  consists of all measurable functions $f$ for which
\[
\|f\|_{p, \alpha}=\left(\int_{{\mathbb C}^n}\left|f(z)\right|^{p}e^
{-\frac{\alpha p}{2}|z|^2}dv(z)\right)^{\frac{1}{p}}<\infty
\]
where $dv$ is  the Lebesgue volume measure on $\mathbb{C}^n$.
The
 classical Fock space  $ F^p_\alpha$ is the subset of $L^{p}_\alpha$, defined to be
  $$
  F^{p}_\alpha=L^{p}_\alpha\cap H({\mathbb C}^n),
  $$
 where $ H({\mathbb C}^n)$ is the family of  all holomorphic function on ${\mathbb C}^n$. And $F^\infty_\alpha$ is the set of those $f\in  H({\mathbb C}^n)$ with
 $$
 \|f\|_{\infty, \alpha}=\sup_{z\in {\mathbb C}^n} |f(z)|e^
{-\frac{\alpha }{2}|z|^2}<\infty.
 $$
 It is clear that, $L^{p}_\alpha$ and  $F^{p}_\alpha$ both are Banach spaces under
 the norm $\|\cdot\|_{p, \alpha}$ if $p\geq 1$.  See \cite{Zh12} for more information.

 A natural generalization involves replacing the Gaussian weight with certain uniformly convex weights \(\varphi\). More precisely, let
 \[\varphi=\varphi(x_{1},x_{2},\ldots,x_{2n})\in C^{2}({\mathbb R}^{2n})\]
  be a real-valued function, and there are positive constants \(m\) and \(M\) such that \(\operatorname{Hess}_{{\mathbb R}}\varphi\), the real Hessian, satisfies
 \begin{equation}\label{a0bwei}
m \textrm{E} \leq\operatorname{Hess}_{{\mathbb R}}\varphi(x)=\left(\frac{\partial^{2}\varphi(x)}{\partial x_{j}\partial x_{k}}\right)_{j,k=1}^{2n}\leq M\textrm{E},
 \end{equation}
where \(\textrm{E}\) is the \(2n\times 2n\) identity matrix. Here, for symmetric matrices \(A\) and \(B\), the notation \(A\leq B\) means that \(B-A\) is positive semidefinite.  It is clear that, if $\varphi(z)=\frac{\alpha}{2}|z|^{2}$, then $\varphi$ satisfies \eqref{a0bwei}. This implies that the theory in \cite{HV22, HV23a, HV23} is applicable in the present setting.
As the argument in \cite{HV23}, a popular example is \(\varphi(z)=|z|^{2}-\frac{1}{2}\log(1+|z|^{2})\), which gives the so-called Fock-Sobolev spaces studied \cite{CZ12}.

Denote by $K (\cdot, \cdot)$  the Bergman kernel of $F^2_\alpha$,  and  $k_{z}$   the normalized reproducing  kernel for $F^2_\alpha$. It is well known that
  $$
K (w, z)=e^{\alpha\langle w, z\rangle} \, \textrm{ and }\, k_z( w)=e^{\alpha\langle w, z\rangle}/e^{\frac{\alpha}{2}|z|^2},\qquad z,w\in\mathbb{C}^n.
$$
Notice that, the set  $$
\textrm{Span} \{k_{z}: z\in {\mathbb C}^n\}
$$
 is dense in $F^{p}_\alpha$ for $p\geq 1$.

 This paper focuses on Hankel operators. Let  $\mathcal S$   be the family of all measurable functions $f$ on ${\mathbb{C}}^n$ satisfying $fk_z\in \cup_{p\ge 1} L^{p}_\alpha$ for each $z\in {\mathbb C}^n$. Given $f \in {\mathcal S}$, the Hankel operator $H_{f}$   is densely defined on $F^{p}_\alpha$ by
\begin{eqnarray*}
H_{f}g=fg-P(fg),
\end{eqnarray*}
where $P$  is the orthogonal projection, defined to be
 $$
 Pf(z)=\int_{{\mathbb C}^n}e^{\alpha\langle z, w\rangle}f(w)e^{-\alpha|w|^2}dv(w).
 $$

Hankel operator is an important model of linear operators in the theory of holomorphic function spaces.  In the past decades, this operator  has been investigated in various contexts on Fock spaces. For example, in  Zhu's book \cite{Zh12},     Hankel operators $H_f$ with   certain kind  of symbols $f$ on  the classical Fock space $F^2_{\alpha}$ are well studied.  For the weight $\varphi$ with \eqref{a0bwei}, Hu and Virtanen  \cite{HV23} obtain the characterizations on $f$ for which the induced Hankel operator $H_f$ is  bounded (or compact) from $F^p_{\varphi}$ to $L^q_{\varphi}$ for $0<p,q<\infty$ recently; they also discuss the Schatten class of Hankel operators on $F^2_{\varphi}$ and obtain the Berger-Coburn phenomenon in \cite{HV22, HV23a}.
However, the theory of absolutely summing Hankel operators remains  undeveloped.\\

  We now introduce the conception of $r$-summing operators.  Suppose $X$ and $Y$ are Banach spaces, and   $1 \leq r < \infty$.  An operator $T: X \to Y$ is called $r$-summing  if it maps weakly $r$-summable sequences in $X$ to absolutely $r$-summable sequences in $Y$, i.e.,  there exists a constant $C \geq 0$ such that for every finite sequence $\{x_k\}_{k=1}^n \subset X$, we have
 \begin{equation}\label{ab-norm}
\left( \sum_{k=1}^n \|T x_k\|_Y^r \right)^{1/r} \leq C \sup_{x^* \in B_{X^*}} \left( \sum_{k=1}^n |x^*(x_k)|^r \right)^{1/r},
\end{equation}
where  $X^*$ denotes the dual space of $X$, and  $B_{X^*} = \{x^* \in X^* : \|x^*\| \leq 1\}$ is the closed unit ball of  $X^*$  (with the weak-$*$ topology).

We shall write $\Pi_r(X, Y )$ for the set of all $r$-summing operators from $X$ to $Y$. The
$r$-summing norm of $T$, denoted by $\pi_r (T : X \to Y )$, is the best constant $C$ for which the inequality
\eqref{ab-norm} always holds. In the following, if no confusion occures, we will abbreviate $\pi_r (T : X \to Y )$ as
$\pi_r (T)$.

The theory of $r$-summing operators has its origins in the pioneering work of Grothendieck during the 1950s. Among his profound contributions is the foundational result that any bounded linear operator mapping from $\ell^1$ into $\ell^2$ is absolutely summing. Interest in this class of operators saw a notable revival in the late 1960s, when their structural properties began to be systematically explored within the framework of functional analysis. A major advancement was made by Pietsch, who formulated a crucial factorization theorem that has since become a cornerstone in the analysis of such operators. Concurrently, Lindenstrauss and Pe{\l}czy\'{n}ski \cite{LP68} clarified the deep connections between summing operators and the geometric structure of Banach spaces.
Subsequent decades have further enriched the theory, as reflected in several influential treatises, including \cite{DJT, Li18, Pe77, Wo91}.

In a recent development, Lef\`{e}vre and Rodr\'{i}guez-Piazza \cite{LR18} established a full characterization of $r$-summing Carleson embeddings on Hardy spaces $H^p$ in the case $p > 1$ and for all $r \geq 1$. He, Jreis, Lef\`{e}vre and Lou  turned attention to the Bergman space setting in \cite{HJLL}. Chen, He and Wang  study absolutely summing   Carleson embeddings on weighted Bergman spaces.  Hu and Wang \cite{HW25} discussed the $r$-summing Toeplitz operators on Fock spaces. Fan  et al. studied the big and little Hankel operators acting between $p$-weighted Bergman and $q$-weighted Lebesgue
spaces  for distinct   exponents $p$ and $q$. It should be noted that in the Bergman space setting, the continuous embedding $A^s_\alpha \hookrightarrow A^t_\alpha$  holds whenever $s > t>0$ and $\alpha > -1$. By contrast, in the setting of Fock spaces, no such continuous embedding exists.\\

For $z\in {\mathbb C}^n$ and $r>0$,   $B(z,r)=\{w\in{\mathbb C}^n: |w-z|<r\}$. Suppose $1 \leq p < \infty$, and suppose $f\in {\mathcal S}$. For $r>0$, define the function $G_{p,r}(f)$ to be
\[
G_{p,r}(f)(z) = \inf_{h\in H(B(z,r))}\left( \frac{1}{|B(z,r)|}\int_{B(z,r)} |f(w)-h(w)|^pdv(w) \right)^{1/p}.
\]
where $H(B(z,r))$ stands for the set of holomorphic functions in  $B(z,r)$.

 As in \cite{HV23}, let the space $\mathrm{IDA}^{s, p}_r$ to be the family of all $f\in L^p_{\text{loc}}$ such that
$$\|f\|_{\mathrm{IDA}^{s, p}_{r}}:=\left(\int_{{\mathbb C}^n}G_{p,r}f(z)^sdv(z)\right)^{1/s}<\infty.$$
By \cite[Corollary 3.10]{HV23}, for $1\leq p<\infty$ and $0<s\leq\infty$, we know that  $\mathrm{IDA}^{s, p}_r$ does not depend on $r$, this means, different values of $r$ give equivalent seminorms $\|G_{p,r}f\|_{L^s}$. We  write
$G_{p}(f)=G_{p, 1}(f)$  and $\mathrm{IDA}^{s, p}=\mathrm{IDA}^{s, p}_{1}$ for simplicity.\\

In this paper, we are going to provide a complete characterization of $f\in {\mathcal S}$ such that,  for all $1 \leq p, r < \infty$, Hankel  operators $H_f$  are $r$-summing from the  Fock space $F^p_\alpha$ to $L^p_\alpha$. For this purpose, denote by $p'$   the conjugate exponent of $p$, and set $\kappa=\kappa (p,r)$ to be
\begin{equation}\label{number-s1}
\kappa =
\begin{cases}
2, & \text{if } 1 \leq p \leq 2 \  \text{and } \ r\geq 1; \\
p', & \text{if } p \geq 2 \  \text{and } \ 1 \leq r \leq p'; \\
r, & \text{if } p \geq 2 \  \text{and } \ p' \leq r \leq p; \\
p, & \text{if } p \geq 2 \  \text{and } \ p \leq r < \infty.
\end{cases}
\end{equation}
The symbol \(A \simeq B\) means \(A \lesssim B \lesssim A\), and \(A \lesssim B\) denotes that there exists a constant \(C > 0\), independent of the relevant variables, such that \(A \leq C B\). Here is the main result.

\begin{theorem}\label{thm:main}  For $1\leq p, r <\infty$ and $f\in {\mathcal S}$, then $H_f\in \Pi_r\left(F^p_\alpha, L^p_\alpha\right)$  if and only if $f \in \mathrm{IDA}^{\kappa,\,p}$.
Furthermore, we have
\[
\pi_r\left(H_f: F_\alpha^p \to L_\alpha^p\right) \simeq\|f\|_{\mathrm{IDA}^{\kappa,\,p}}.
\]
\end{theorem}

The paper is organized as follows. In Section 2, we collect some preliminary knowledge that will be used in later sections.   Section 3 is devoted to proving the main result,  about $r$-summing Hankel operators on Fock spaces.  Our proof will be carried out by considering several different cases  based on the values of $p$ and $r$. In Section 4, we study the  Berger-Coburn phenomenon for absolutely summing-Hankel operators on Fock spaces.

\section{Some preliminaries }\label{sect2}

In this section, we share  some fundamental  preliminaries used in the proof of the main theorem. First, we give some notation.

 Suppose $\{a_k\}_{k=1}^\infty$ is a sequence in  ${\mathbb C}^n$. For  $r>0$, we call
$\{a_k\}$  is  an $r$-lattice if the balls $\left\{B(a_k, r)\right\}_{k=1}^{\infty}$ cover
 ${\mathbb C}^n$ and
 $\left\{B\left(a_k, br\right)\right\}_{k=1}^{\infty}$ are pairwise disjoin for some $b>0$.
Via standard treatment, for $c\geq 1$,  there exists some
integer $N$ such that each $z\in {\mathbb C}^n$ can be in at most $N$ balls of $\{B(a_k, cr)\}$. Equivalently,
\begin{equation}\label{number}
1\leq \sum\limits_{k=1}^{\infty}\chi_{B(a_k,cr)}(z) \leq N.
\end{equation}

 Now we turn to introduce some vector-valued sequence spaces, which can be found in \cite{DJT}. Suppose   $1 \leq p < \infty$ and $X$ is a Banach space. A vector sequence $\{x_k\}_{k=1}^\infty$ in $X$ is called strongly $p$-summable (denoted by $\{x_k\}\in \ell_{p}^{\text{strong}}(X)$) if the corresponding scalar sequence $\{\|x_k\|_X\}_{k=1}^\infty$ belongs to $l^p$, and a natural norm is given by
\[
\left \|\{x_k\}\right \|_{p}^{\text{strong}} := \left( \sum_{k=1}^\infty \|x_k\|_X^p \right)^{1/p}.
\]

Strong $p$-summability makes reference to the strong (or norm) topology on $X$. What about the natural analogue for the weak topology? A vector sequence $\{x_k\}_{k=1}^\infty$ in $X$ is called weakly $p$-summable (denoted by $\{x_k\}\in  \ell_p^{\text{weak}}(X)$) if the scalar sequences $\{\langle x^*, x_k\rangle\} \in l^p$ for every $x^* \in X^*$. Of course,  $\ell_p^{\text{weak}}(X)$ is a Banach space with the norm
\[
\|\{x_k\}\|_p^{\text{weak}} := \sup \left\{ \left( \sum_{k=1}^\infty  |\langle x^*, x_k \rangle|^p \right)^{1/p} : x^* \in B_{X^*} \right\}.
\]

For  Banach spaces \(X\) and \(Y\), let \({\mathcal{B}}(X,Y)\) be the collection of all bounded linear operators from \(X\) to \(Y\). For \(T\in \mathcal{B}(X,Y)\), the correspondence
\[\hat{T}:\{x_{k}\}_{k}\mapsto\{Tx_{k}\}_{k}\]
always induces a bounded linear operator from  $\ell_{p}^{\mathrm{weak}}(X)$ to $\ell_{p}^{\mathrm{weak}}(Y)$, as well as a bounded linear operator from  $\ell_{p}^{\mathrm{strong}}(X)$ to $\ell_{p}^{\mathrm{strong}}(Y)$. Moreover,
$$
\|T\|_{X \to Y}=\|\hat{T}\|_{\ell_{p}^{\mathrm{weak}}(X)\to \ell_{p}^{\mathrm{weak}}(Y)}=\|\hat{T}\|_{\ell_{p}^{\mathrm{strong}}(X)\to \ell_{p}^{\mathrm{strong}}(Y)},
$$
which can be found in \cite[Page 34]{DJT}.

It is natural to ask whether $\hat{T}$ is bounded from $\ell_{p}^{\mathrm{weak}}(X)$ to $\ell_{p}^{\mathrm{strong}}(Y)$.  From \cite[Proposition 2.1]{DJT}, we know
\[
T\in \Pi_{p}(X, Y) \ \ {\rm if\ and\ only\ if}\ \  \hat{T}\in \mathcal{B}(\ell_{p}^{\mathrm{weak}}(X)\to \ell_{p}^{\mathrm{strong}}(Y))
\]
 with
\[
 \pi_{p}(T: X \to Y)\simeq \|\hat{T}\|_{\ell_{p}^{\mathrm{weak}}(X)\to \ell_{p}^{\mathrm{strong}}(Y)}.
\]

More generally, for $1 \leq p, q < \infty$, the operator $T: X \to Y$ is called $(q, p)$-summing if the induced operator
$\hat{T}$ is bounded from $\ell_{p}^{\text{weak}}(X)$ to $\ell_{q}^{\text{strong}}(Y)$. Such operators form a vector space, denoted by
\[
\Pi_{q,p}(X, Y),
\]
which is a Banach space under the norm
\[
\pi_{q,p}(T: X \to Y) := \|\hat{T}\|_{\ell_{p}^{\text{weak}}(X) \to \ell_{q}^{\text{strong}}(Y)}.
\]
A standard characterization shows that, a bounded linear operator $T$ belongs to $\Pi_{q,p}(X, Y)$ if and only if there exists a constant $K \geq 0$ such that
\[
\left( \sum_{k=1}^n \|T x_k\|_Y^q \right)^{1/q} \leq K \cdot \sup \left\{ \left( \sum_{k=1}^n |\langle x^*, x_k \rangle|^p \right)^{1/p} : x^* \in B_{X^*} \right\}
\]
holds for every finite collection $\{x_1, \dots, x_n\} \subset X$. Moreover, $\pi_{q,p}(u)$ equals to the smallest such constant $K$.
It is trivial that,
 $$
\Pi_{p,p}(X, Y)=\Pi_{p}(X, Y).
 $$

From \cite{DJT}, we also know

(1) For $1 \leq p < \infty$, the Lebesgue space $L^p(\Omega, \mu)$ has cotype $\max\{p, 2\}$. See Page 219;

(2) For $1 \leq p<q < \infty$,
 $
 \Pi_p(X, Y)\subset \Pi_q(X, Y)
 $ with $\pi_q(T)\leq \pi_p(T)$ for $T\in  \Pi_p(X, Y)$. See Page 39;

 (3) For Banach spaces $X$ and $Y$, if $X$ and $Y$ both have cotype 2,
 \[
 \Pi_r(X, Y)= \Pi_1(X, Y)
 \]
with the equivalence $\pi_r(T) \simeq \pi_1(T)$ for $1 \leq r < \infty$. See Corollary 11.16;

(4) For Banach spaces $X$ and $Y$,  if $Y$ has cotype $q\in [2, \infty)$, then
\[
\Pi_r(X,Y) \subset \Pi_{q,2}(X,Y)
\]
for  all $1 \leq r < \infty$. In particular, if $Y$ has cotype 2, then for any Banach space $X$ and all $2 < r < \infty$, $\Pi_r(X,Y) = \Pi_2(X,Y)$. See Theorem 11.13.\\

Here are some  lemmas for our later use. Lemma \ref{lemma1} and Lemma \ref{lemma3}  can be found in \cite{DJT}. Lemma \ref{lemma2} comes from
\cite[Lemma 2.4]{HL14}.
 These  results collectively establish  the framework needed to prove the main theorem relating the $r$-summing norm of $H_f$ to the $\mathrm{IDA}^{s, p}$-norm of $f$.

\begin{lemma}[Khintchine's Inequality] \label{lemma1}
Let $\{\gamma_j(t)\}_{j=1}^\infty$ be a sequence of Rademacher functions on $[0,1]$, and let $0 < p < \infty$. Then, for any finite sequence of complex numbers $\{c_j\}_{j=1}^n$, we have
\[
\left( \int_0^1 \left| \sum_{j=1}^n c_j \gamma_j(t) \right|^p dt \right)^{1/p} \simeq \left( \sum_{j=1}^n |c_j|^2 \right)^{1/2}.
\]
\end{lemma}

\begin{lemma}  \label{lemma2}
Let $\{c_j\}_{j=1}^\infty \in l^p$ with $0<p<\infty$, and let  $\{z_j\}_{j=1}^\infty \subset \C^n$ be a $r$-lattice. Define
\[
g(z) = \sum_{j=1}^\infty c_j  k_{z_j}(z), \qquad z\in {\mathbb C}^n.
\]
Then $g\in F^p_\alpha$ with
\[
 \|g\|_{p, \alpha}\lesssim \left(\sum_{j=1}^\infty |c_j|^p\right)^{1/p}.
\]
\end{lemma}

\begin{lemma}[Pietsch Domination Theorem]  \label{lemma3}
Suppose $1 \leq r< \infty$, and suppose $X, Y$ both are  Banach spaces.  An operator $T: X \to Y$ is $r$-summing if and only if, there exists a regular Borel probability measure $\sigma$ on the closed unit ball $B_{X^*}$ such that for all $x \in X$, we have
\[
\|T(x)\|_Y \leq \pi_r(T) \left( \int_{B_{X^*}} |x^*(x)|^r d\sigma(x^*) \right)^{1/r}.
\]
\end{lemma}

Suppose $1\leq p<\infty$, and suppose $\mu$ is  a positive Borel measure on ${\mathbb C}^n$. The  space
$L^{p}_\alpha(d\mu)$ is the family of all $\mu$-measurable functions $f$ on ${\mathbb C}^n$ for which
$$
\|f\|_{L^p_\alpha(d\mu)}=\left(\int_{{\mathbb C}^n}\left|f(z)\right|^{p}e^
{-\frac{\alpha p}{2}|z|^2}d\mu(z)\right)^{\frac{1}{p}}<\infty.
$$
In \cite{HL14}, we discussed some characterizations on  $\mu$ such that the embedding operator $\operatorname{Id}$ is bounded or compact from $F^p_{\alpha}$ to $L^p_\alpha(d\mu)$ in terms of the averaged function $\widehat{\mu}$ or the $t$-Berezin transform $\widetilde{\mu}_t$ ($t>0$), where
\[
\widehat{\mu}(z) = \frac{\mu(B(z,1))}{|B(z,1)|} \quad \text{and} \quad  \widetilde{\mu}_t(z) = \int_{{\mathbb C}^n}|k_z(w)|^te^{-\frac{\alpha t}{2}|w|^2}d\mu(w)
\]
for $z \in \mathbb{C}^n$.
Now we are in the position to state the property of absolutely summing embedding operators from $F^p_{\alpha}$ to $L^p_\alpha(d\mu)$.

\begin{proposition} [Absolutely Summing Carleson Embedding] \label{embedding} Suppose  $\mu$ is a positive Borel measure on $\mathbb{C}^n$. Let $1\leq p<\infty$, and let $\kappa$ be as \eqref{number-s1}.
Then the embedding operator $\operatorname{Id}\in \Pi_r\left(F^p_{\alpha}, L^p_\alpha(d\mu)\right)$  if and only if $\widehat{\mu} \in L^{\kappa/p}(dv)$ $($or $\widetilde{\mu}_t \in L^{\kappa/p}(dv)$ for any $t>0)$.
 Moreover,
\[
\pi_r\left(\operatorname{Id}:F^p_{\alpha} \to L^p_\alpha(d\mu)\right) \simeq \|\widehat{\mu}\|^{1/p}_{L^{\kappa/p}}  \simeq \|\widetilde{\mu}_t\|^{1/p}_{L^{\kappa/p}} .
\]
\end{proposition}

\begin{proof}
The $L^{\kappa/p}$-norm equivalence between $\widehat{\mu}$ and $\widetilde{\mu}_t$ follows from \cite[Theorem 2.8]{HL14}.
Similar to the proof in \cite{CHW} or \cite{HJLL}, we can  prove the result in the case of \(1 < p <\infty\).
We now deal with the case $p=1$.
First, we assume $\widehat{\mu}\in L^2({\mathbb C}^n)$. Consider the commutative diagram
$$
\xymatrix{
 {F}_{\alpha}^{1}\ar[d]_{\operatorname{Id}} \ar[r]^{\operatorname{Id}} &   L^1_\alpha(d\mu)       \\
 {F}_{\alpha}^{2} \ar[ur]_{\operatorname{Id}}                     }
 $$
By \cite[Theorem 2.8]{HL14}, we know that $\mu$ is a  (2,1) Fock-Carleson measure, then
\[
\|\operatorname{Id}\|_{F^2_\alpha \to L^1_\alpha(d\mu)} \lesssim\norm{\widehat{\mu}}_{L^2}.
\]
Since
\[
\pi_1(\operatorname{Id}: F^1_\alpha \to F^2_\alpha)  \lesssim \|\operatorname{Id}\|_{F^1_\alpha \to F^2_\alpha} <\infty,
\]
we get
\[
\pi_1(\operatorname{Id}: F^1_\alpha \to L^1_\alpha(d\mu))  \lesssim \pi_1(\operatorname{Id}: F^1_\alpha \to F^2_\alpha) \cdot \|\operatorname{Id}\|_{F^2_\alpha \to L^1_\alpha(d\mu)}  \lesssim \|\widehat{\mu}\|_{L^2}.
\]

Conversely, suppose $\pi_1(\operatorname{Id}: F^1_\alpha \to L^1_\alpha(d\mu)) < \infty$. Notice that,
$(F^1_\alpha)^*=F^\infty_\alpha$ under the pairing $$\langle f, g\rangle=\int_{{\mathbb C}^n}f(w)\overline{g(w)}e^{-\alpha|w|^2}dv(w).$$
By Pietsch's Domination Theorem, there exists a probability measure $\sigma$ on $B_{\left(F^1_\alpha\right)^*}$ such that for all $g \in F^1_\alpha$, we obtain
\[
\norm{\operatorname{Id}(g)}_{L^1_\alpha(d\mu)} \leq \pi_1\left(\operatorname{Id}: F^1_\alpha \to L^1_\alpha(d\mu) \right) \int_{B_{\left(F^1_\alpha\right)^*}} \abs{\langle h, g\rangle} d\sigma
(h).
\]
  Let $\{c_j\}_{j=1}^\infty$ be any complex sequence in $l^{2}$.  For fixing positive integer $m$, define the random function
 \begin{equation}\label{textfunct}
g_t(w) = \sum_{j=1}^m c_j \gamma_j(t)k_{z_j}(w), \quad w\in {\mathbb C}^n,
\end{equation}
where   $\{\gamma_j(t)\}$ is a sequence of Rademacher functions, and $\{z_j\}$ is a $\delta$-lattice on ${\mathbb C}^n$, $\delta>0$ is small enough. Then $g_t\in F^1_\alpha$ and
$$
\|g_t\|_{F^1_\alpha}\lesssim \sum_{j=1}^m |c_j|.
$$
By the property of the reproducing kernel for $F^\infty_\alpha$,  we have
\[
\int_{B_{\left(F^1_\alpha\right)^*}} \left|\langle h, g_t\rangle \right|d\sigma(h) = \int_{B_{F^\infty_\alpha}} \left|\sum_{j=1}^m c_j \gamma_j(t) h(z_j)e^{-\frac{\alpha}{2}|z_j|^2}\right|d\sigma
(h).
\]
Using Khinchin's inequality, we obtain
\begin{align*}
 \int_0^1 \left( \int_{B_{F^\infty_\alpha}} \left|\langle h, g_t\rangle\right| d\sigma(h) \right) dt
&=
\int_{B_{F^\infty_\alpha}} \left( \int_0^1 \left|\sum_{j=1}^m c_j \gamma_j(t) h(z_j)e^{-\frac{\alpha}{2}|z_j|^2}\right| dt \right) d\sigma(h) \\
&
\leq \int_{B_{F^\infty_\alpha}} \left( \sum_{j=1}^m \left|c_j h(z_j)e^{-\frac{\alpha}{2}|z_j|^2}\right|^2 \right)^{1/2} d\sigma(h) \\
&
\lesssim \int_{B_{F^\infty_\alpha}} \norm{h}_{\infty, \alpha} \left( \sum_{j=1}^m \left|c_j\right|^2 \right)^{1/2} d\sigma(h) \\
&
\leq \norm{\{c_j\}}_{l^2}.
\end{align*}
On the other hand, Khinchin's inequality and \eqref{number} show
\begin{align*}
 & \ \ \ \   \int_0^1 \norm{\operatorname{Id}(g_t)}_{L^1_\alpha(\mu
)} dt\\
&
\simeq \int_{{\mathbb C}^n} \left( \sum_{j=1}^m \left|c_j k_{z_j}(w)\right|^2 \right)^{1/2} e^{-\frac{\alpha}{2}|w|^2} d\mu(w) \\
&
\geq \frac{1}{N}\sum_{k=1}^\infty \int_{B(z_k, \delta)} \left( \sum_{j=1}^m \left|c_j\right|^2 \left|k_{z_j}(w)\right|^2 \right)^{1/2} e^{-\frac{\alpha}{2}|w|^2} d\mu(w) \\
&
\geq \frac{1}{N}\sum_{j=1}^m \left|c_j\right| \int_{B(z_j, \delta)} \left|e^{\alpha\langle w, z_j\rangle}\right|e^{-\frac{\alpha}{2}|z_j|^2} e^{-\frac{\alpha}{2}|w|^2} d\mu(w) \\
&
\geq \frac{1}{N}\sum_{j=1}^m \left|c_j\right| \mu(B(z_j,  \delta))\cdot \inf_{w \in B(z_j,  \delta)} e^{-\frac{\alpha}{2}| w- z_j|^2} \\
&
\geq \frac{C}{N}\sum_{j=1}^m \left|c_j\right| \whmu
_{\delta}(z_j).
\end{align*}
Hence,
\[
\sum_{j=1}^m \left|c_j\right| \whmu_{\delta}(z_j) \lesssim \pi_1\left(\operatorname{Id}: F^1_\alpha \to L^1_\alpha(d\mu) \right) \norm{\{c_k\}
}_{l^2}.
\]
Letting $m\rightarrow\infty$, by the duality, we obtain
\[
\left(\sum_{j=1}^\infty  \whmu_{\delta}(z_j)^2\right)^{1/2} \lesssim \pi_1\left(\operatorname{Id}: F^1_\alpha \to L^1_\alpha(d\mu) \right).
\]
Using \cite[Theorem 2.8]{HL14} again, we have
 \[
\left( \int_{{\mathbb C}^n} \left|\whmu(w)\right|^2 dv(w) \right)^{1/2} \simeq \left(\sum_{j=1}^\infty  \whmu_{\delta}(z_j)^2\right)^{1/2}  \lesssim \pi_1\left(\operatorname{Id}: F^1_\alpha \to L^1_\alpha(d\mu) \right).
\]
Therefore, $$
\pi_1\left(\operatorname{Id}: F^1_\alpha \to L^1_\alpha(d\mu) \right)\simeq\|\widehat{\mu}\|_{L^2}.
$$
Notice that, $F^1_\alpha$   and  $L^1_\alpha(d\mu)$  both have cotype 2, so
  $$\pi_r\left(\operatorname{Id}: F^1_\alpha \to L^1_\alpha(d\mu) \right) \simeq \pi_1\left(\operatorname{Id}: F^1_\alpha \to L^1_\alpha(d\mu) \right)\simeq\|\widehat{\mu}\|_{L^2}$$ for $1 < r < \infty$.
This completes the proof.
\end{proof}

\section{Proof of Theorem  \ref{thm:main} }\label{sect3}

The proof of Theorem \ref{thm:main} will be contained in the following three theorems in this section.


\begin{theorem}\label{thm:main-1}
 Suppose $f\in {\mathcal S}$.  Let $1\leq p<\infty$, and let $\kappa $ be  as  \eqref{number-s1}. If
  $f\in \mathrm{IDA}^{\kappa,\,p}$,  then the Hankel operator $H_f: F^p_\alpha \to L^p_\alpha$ is $r$-summing. Moreover, we have
\[
\pi_r(H_f: F^p_\alpha \to L^p_\alpha) \lesssim  \|f\|_{\mathrm{IDA}^{\kappa,\,p}}.
\]
\end{theorem}

\begin{proof}
First, we  assume $f\in \mathrm{IDA}^{\kappa,\,p}$, then $G_{p, \tau}(f)\in L^{\infty}(dv)$ and $G_{p, \tau}(f)\in L^{\kappa}(dv)$ for all $\tau>0$. The proof of  \cite[Theorem 3.8]{HV23} imples that  there exist $f_1$ and $f_2$ such that
$f=f_1+f_2$ satisfying
\begin{equation}\label{basic-2.1}
 \|M_{p,1}\left(\overline{\partial}f_1\right)\|
 _{L^\kappa}+\|M_{p,1}\left(f_2\right)\|_{L^\kappa}\lesssim\|G_{p, \tau}(f)\|_{L^\kappa}
\end{equation}
for some $\tau>0$,  where $$M_{p,1}\left(f\right)(z)=\left(\frac{1}{|B(z,1)|}\int_{B(z,1)}|f|^pdv\right)^{1/p}.$$ By the result on Page 2062 in Hu and Virtanen \cite{HV23}, we know $H_f$ is bounded from $F^p_\alpha$ to $L^p_\alpha$, and   for any $g\in F^p_\alpha$,
\[
\|H_{f_1}(g)\|_{p, \alpha} \lesssim\left\| g\left|\overline{\partial}f_1\right|\right\|_{p, \alpha}  \  \text{and }  \
\|H_{f_2}(g)\|_{p, \alpha}  \lesssim\|f_2 g\|_{p, \alpha}.
\]
It follows that
\[
\pi_r(H_{f_1})\lesssim \pi_r\left(M_{\left|\overline{\partial}f_1\right|}\right) \  \text{and }  \  \pi_r\left(H_{f_2}\right)\lesssim \pi_r\left(M_{f_2}\right),
\]
where $M_f$ denotes the pointwise multiplier, that is $M_{f} g=fg$.
Set $$d\mu_1=\left|\overline{\partial}f_1\right|^pdv  \  \text{and }  \  d\mu_2=\left|f_2\right|^pdv.$$
By Proposition \ref{embedding}, we have
\[
  \pi_r\left(M_{\left|\overline{\partial}f_1\right|}\right)=  \pi_r\left(\operatorname{Id}: F^p_\alpha\rightarrow L^p_\alpha(d\mu_1)\right) \simeq  \|\widehat{\mu}_1\|^{1/p}_{L^{\kappa/p}} \simeq  \|M_{p,1}\left(\overline{\partial}f_1\right)\|_{L^\kappa}
    \]
  and
   \[
  \pi_r\left(M_{f_2}\right)=  \pi_r\left(\operatorname{Id}: F^p_\alpha\rightarrow L^p_\alpha(d\mu_2)\right) \simeq  \|\widehat{\mu}_2\|^{1/p}_{L^{\kappa/p}}\simeq  \|M_{p,1}\left(f_2\right)\|_{L^\kappa}.
    \]
     \eqref{basic-2.1} yields
    \[
  \pi_r\left(H_f\right)\leq\pi_r\left(H_{f_1}\right)+\pi_r\left(H_{f_2}\right) \lesssim  \|G_{p, \tau}\left(f\right)\|_{L^\kappa} \simeq \|f\|_{\mathrm{IDA}^{\kappa,\,p}}.
    \]
   The proof is complete.
\end{proof}

\begin{theorem}\label{thm:main-2}
Suppose $1\leq p<2$, and suppose  $H_f: F^p_\alpha \to L^p_\alpha$ is $r$-summing for some $r\geq 1$. Then $f\in  \mathrm{IDA}^{2,\,p}$. Moreover, we have
\[
\|f\|_{\mathrm{IDA}^{2,  p}}\lesssim\pi_r(H_f: F^p_\alpha \to L^p_\alpha).
\]
\end{theorem}

\begin{proof}
Suppose $H_f$ is $r$-summing. Since $F^p_\alpha$ and $L^p_\alpha$ both
have cotype  $2$ when $1\leq p< 2$, \cite[Corollary 11.16]{DJT} shows
$$
H_f\in \Pi_r(F^p_\alpha, L^p_\alpha) = \Pi_p(F^p_\alpha, L^p_\alpha)=\Pi_1(F^p_\alpha, L^p_\alpha)
$$
 with the equivalence
  $$
 \pi_r(H_f) \simeq\pi_p(H_f)\simeq\pi_1(H_f).
  $$
 Let $\{c_j\}_{j=1}^\infty$ be any complex sequence in $l^{\frac{2p}{2-p}}$.  Suppose $\{\gamma_j(t)\}$ is a sequence of Rademacher functions, and $\{z_j\}$ is a $\delta$-lattice on ${\mathbb C}^n$, where $\delta>0$ is small enough.  For fixing positive integer $m$, define the random function as in \eqref{textfunct}.
 Then $g_t\in F^p_\alpha$ and
$$
\|g_t\|_{p, \alpha}\lesssim \left(\sum_{j=1}^m |c_j|^p\right)^{1/p}.
$$
So
$$
\int_0^1 \|H_f(g_t)\|_{p, \alpha}^p \, dt\leq \|H_f\|_{F^p_\alpha\rightarrow L^p_\alpha}^p  \|g_t\|^p_{p, \alpha}<\infty.
$$
Using Fubini's theorem and Khintchine's inequality, we have
\begin{align*}
&\int_0^1 \|H_f(g_t)\|_{p, \alpha}^p \, dt\\
&= \int_{{\mathbb C}^n} \left( \int_0^1 \left| \sum_{j=1}^m c_j \gamma_j(t) H_f(k_{z_j})(w) \right|^p dt \right) e^{-\frac{\alpha p}{2}|w|^2} dv(w) \\
&\simeq \int_{{\mathbb C}^n}  \left( \sum_{j=1}^m |c_j|^2 |H_f(k_{z_j})(w)|^2 \right)^{p/2}  e^{-\frac{\alpha p}{2}|w|^2} dv(w) \\
&\geq \sum_{j=1}^m |c_j|^p \int_{B(z_j,2\delta)} |H_f(k_{z_j})(w)|^p  e^{-\frac{\alpha p}{2}|w|^2} dv(w).
\end{align*}
It is well known that
   $\inf\limits_{w \in B(z, 2\delta)} |k_{ z}(w)|e^{-\frac{\alpha}{2}|w|^2} >0$,
  which indicates
   $$\frac{P(f k_{z_j})}{k_{ z_j} }\in H(B(z_j, 2\delta))$$
    for all $j$. This gives
  \begin{align*}
G_{p, 2\delta}(f)(z_j)^p \lesssim\int_{B(z_j,2\delta)} |H_f(k_{z_j})(w)|^p  e^{-\frac{\alpha p}{2}|w|^2} dv(w).
\end{align*}
Thus,
 \begin{align*}
\sum_{j=1}^m |c_j|^pG_{p, 2\delta}(f)(z_j)^p \lesssim\int_0^1 \|H_f(g_t)\|_{p, \alpha}^p \, dt<\infty.
\end{align*}

On the other hand, by the Pietsch domination theorem for $r$-summing operators, there exists a probability measure $\sigma$ on the unit ball of $(F^p_\alpha)^{*}=F^{p'}_\alpha$ (denoted by $B_{F^{p'}_\alpha}$) such that
\begin{align*}
\|H_f(g_t)\|_{p, \alpha}^p &\leq \pi_p(H_f)^p \cdot\left( \int_{B_{F^{p'}_\alpha}} \left| \langle h, g_t \rangle \right|^p d\sigma(h) \right)\\
&\simeq \pi_r(H_f)^p\cdot\int_{B_{F^{p'}_\alpha}} \left| \sum_{j=1}^{m} c_j \gamma_j(t) h(z_j)e^{-\frac{\alpha}{2}|z_j|^2} \right|^p  d\sigma(h).
\end{align*}
The last  equality is due to the property of the reproducing kernel for $F^{p'}_\alpha$. Hence,
\begin{align*}
&\int_0^1 \| H_f (g_t) \|_{p, \alpha}^p dt\\
&\leq \pi_r(H_f)^p\cdot\int_{B_{F^{p'}_\alpha}} d\sigma(h) \int_0^1  \left| \sum_{j=1}^{m} c_j \gamma_j(t) h(z_j)e^{-\frac{\alpha}{2}|z_j|^2} \right|^p dt \\
&\simeq \pi_r(H_f)^p\cdot\int_{B_{F^{p'}_\alpha}}  \left( \sum_{j=1}^m |c_j|^2 |h(z_j)|^2e^{-\alpha|z_j|^2} \right)^{\frac{p}{2}} d\sigma(h).
\end{align*}
 Notice that, $\frac{p'}{2}=\frac{p}{2(p-1)}>1$ and $\frac{p}{2-p}$ is the conjugate exponent of $\frac{p'}{2}$. By H\"{o}lder's inequality, we obtain
\begin{align*}
&\left(\sum_{j=1}^m |c_j|^2 |h(z_j)|^2e^{-\alpha|z_j|^2}  \right)^{\frac{p}{2}} \\
&\leq \left[\sum_{j=1}^m  \left(|h(z_j)|^2e^{-\alpha|z_j|^2} \right)^{\frac{p'}{2}}\right]^{\frac{2}{p'}\cdot\frac{p}{2}}  \left(\sum_{j=1}^m |c_j|^{\frac{2p}{2-p}}\right)^{\frac{2-p}{p}\cdot\frac{p}{2}}\\
&=\left[\sum_{j=1}^m  \left|h(z_j)e^{-\frac{\alpha}{2}|z_j|^2}\right|^{p'}\right]^{\frac{p}{p'}}  \left(\sum_{j=1}^m |c_j|^{\frac{2p}{2-p}}\right)^{\frac{2-p}{2}}
\end{align*}
Since $\{z_j\}$ is a $\delta$-lattice on ${\mathbb C}^n$, we get
\begin{align*}
\sum_{j=1}^m  \left|h(z_j)e^{-\frac{\alpha}{2}|z_j|^2}\right|^{p'}&\lesssim\sum_{j=1}^m \int_{B(z_j, \delta)}\left|h(w)e^{-\frac{\alpha}{2}|w|^2}\right|^{p'}dv(w)\\
&\leq N\int_{{\mathbb C}^n}\left|h(w)e^{-\frac{\alpha}{2}|w|^2}\right|^{p'}dv(w).
\end{align*}
Therefore,
\begin{align*}
\int_0^1 \| H_f (g_t) \|_{p, \alpha}^p dt
&\lesssim \pi_r(H_f)^p\cdot\int_{B_{F^{p'}_\alpha}} \|h\|^p_{p^{'},\,\alpha} \cdot \left(\sum_{j=1}^m |c_j|^{\frac{2p}{2-p}}\right)^{\frac{2-p}{2}} d\sigma(h)\\
& \leq \pi_r(H_f)^p\cdot\| \{c_j\}_{j=1}^m \|_{l^{\frac{2p}{2-p}}}^p.
\end{align*}
Combining the above statistics, we conclude
 \begin{align*}
\sum_{j=1}^m |c_j|^pG_{p, 2\delta}(f)(z_j)^p \lesssim\int_0^1 \| H_f (g_t) \|_{p, \alpha}^p dt \lesssim\pi_r(H_f)^p\| \{c_j\}_{j=1}^\infty \|_{l^{\frac{2p}{2-p}}}^p.
\end{align*}
Note that,  $\frac{2}{p}$ is the   conjugate exponent of $\frac{2}{2-p}$. Taking $m\rightarrow\infty$, by the duality, we obtain
\[ \left(\sum_{j=1}^\infty G_{p, 2\delta}(f)(z_j)^2 \right)^{1/2}\lesssim \pi_r(H_f). \]
By the definition, we have
 \begin{equation}\label{norm02}
G_{p, \delta} (f)(z) \lesssim G_{p, 2\delta} (f)(w), \  \ \ z\in B(w, \delta).
 \end{equation}
We have
\[
 \int_{{\mathbb C}^n} G_{p, \delta}(f)^2dv\leq \sum_{j=1}^\infty \int_{B(z_j, \delta)}G_{p, \delta}(f)^2dv \lesssim \sum_{j=1}^\infty G_{p, 2\delta}(f)(z_j)^2.
\]
This shows
\[
 \int_{{\mathbb C}^n} G_p(f)^2 dv\simeq\int_{{\mathbb C}^n} G_{p, \delta}(f)^2dv \lesssim \sum_{j=1}^\infty G_{p, 2\delta}(f)(z_j)^2\lesssim \pi_r(H_f)^2.
\]
This completes the proof.
\end{proof}

\begin{theorem}\label{thm:main-3}
Let $p\geq 2$, and let $\kappa$ be as in (\ref{number-s1}).
If   $H_f: F^p_\alpha \to L^p_\alpha$ is $r$-summing, then $f \in \mathrm{IDA}^{\kappa,\,p}$. Moreover, we have
\[
\|f\|_{\mathrm{IDA}^{\kappa,\,p}}\lesssim\pi_r(H_f: F^p_\alpha \to L^p_\alpha).
\]
\end{theorem}

\begin{proof} We separate the proof into three cases.

{\bf Case 1}:   $1 \leq r \leq p'$.

Suppose $H_f\in \Pi_r(F^p_\alpha, L^p_\alpha)$, then    $\pi_{p'}(H_f)\leq \pi_r(H_f)<\infty$, since
$$
\Pi_{r}(F^p_\alpha, L^p_\alpha)\subseteq \Pi_{p'}(F^p_\alpha, L^p_\alpha) \  \text{ for } \ 1 \leq r \leq p'.
$$
For any $\delta$-lattice  $\{z_j\}$ and any fixed
integer $m>0$, we know $\left\{k_{z_j}\right\}_{j=1}^{m}$ is a finite sequence in $F^{p}_\alpha$. By the definition of $p'$-summing operator,  we have
\begin{align*}
\sum\limits_{j=1}^{m}\|H_f(k_{z_j})\|_{p, \alpha}^{p'} &\leq \pi_{p'}(H_f)^{p'}
\cdot\sup_{h\in B_{F^{p'}_\alpha}} \sum\limits_{j=1}^{m}\left| \langle h, k_{z_j} \rangle \right|^{p'} \\
&\lesssim\pi_{r}(H_f)^{p'} \cdot \sup_{h\in B_{F^{p'}_\alpha}} \sum\limits_{j=1}^{m}\left|h(z_j)e^{-\frac{\alpha}{2}|z_j|^2}\right|^{p'},
\end{align*}
the last inequality comes from the the property of the reproducing kernel for  $F^{p'}_\alpha$. Notice that,
\begin{align*}
\sum\limits_{j=1}^{m}\left|h(z_j)e^{-\frac{\alpha}{2}|z_j|^2}\right|^{p'}&\lesssim\sum\limits_{j=1}^{m}\frac{1}{|B(z_j, \delta)|}\int_{B(z_j, \delta/2)}\left|h(w)e^{-\frac{\alpha }{2}|w|^2}\right|^{p'}dv(w)\\
&\leq N\int_{{\mathbb C}^n}\left|h(w)e^{-\frac{\alpha }{2}|w|^2}\right|^{p'}dv(w)\leq N
\end{align*}
for $h\in B_{F^{p'}_\alpha}$. This implies
\begin{align*}
\sum\limits_{j=1}^{m}\|H_f(k_{z_j})\|_{p, \alpha}^{p'} \lesssim \pi_{r}(H_f)^{p'}.
\end{align*}
On the other hand, similarly to the proof in Theorem \ref{thm:main-2}, we get
\begin{align*}
\|H_f(k_{z_j})\|_{p, \alpha}^{p'}
&\geq \left(\int_{B(z_j, 2\delta)} \left|k_{z_j}f(w)-P(k_{z_j}f)(w)\right|^p  e^{-\frac{\alpha p}{2}|w|^2} dv(w)\right)^{p'/p}\\
&\simeq\left(\frac{1}{|B(z_j, 2\delta)|}\int_{B(z_j,2\delta)}\left |f(w)-\frac{P(k_{z_j}f)(w)}{k_{z_j}(w)}\right|^p  dv(w)\right)^{p'/p}\\
&\geq G_{p,2\delta}(f)(z_j)^{p'}.
\end{align*}
 These show
  \begin{align*}
\sum\limits_{j=1}^{m}G_{p, 2\delta}(f)(z_j)^{p'} \lesssim \sum\limits_{j=1}^{m}\|H_f(k_{z_j})\|_{p, \alpha}^{p'}
\lesssim \pi_{r}(H_f)^{p'} .
\end{align*}
Letting $m\rightarrow\infty$, \eqref{norm02} gives
 \begin{align*}
\int_{{\mathbb C}^n}G_p(f)^{p'}dv&\leq\sum\limits_{j=1}^{\infty}\int_{B(z_j,\delta)}G_{p, \delta}(f)^{p'}dv\\
&\lesssim \sum\limits_{j=1}^{\infty}G_{p, 2\delta}(f)(z_j)^{p'}\\
&\lesssim \pi_{r}(H_f)^{p'}<\infty.
\end{align*}

{\bf Case 2}:  $p' \leq r \leq p$.

Suppose $H_f: F^p_\alpha \to L^p_\alpha$ is $r$-summing. Similarly to the proof in Case 1, we obtain
\begin{align*}
\int_{{\mathbb C}^n}G_p(f)^{r}dv\lesssim \sum\limits_{j=1}^{\infty}G_{p, 2\delta}(f)(z_j)^{r}
\lesssim \pi_{r}(H_f)^r<\infty,
\end{align*}
where $\{z_j\}$ is a $\delta$-lattice on ${\mathbb C}^n$.

{\bf Case 3}:  $p \leq r <\infty$.

Suppose  $H_f\in \Pi_r(F^p_\alpha, L^p_\alpha)$.  Since the cotype of $L^p_\alpha$
is $p \geq 2$, it
follows from \cite[Theorem 11.13]{DJT} that
 $H_f: F^p_\alpha \to L^p_\alpha$
is $(p, 2)$-summing, and
 $$
\pi_{p, 2}(H_f)\lesssim \pi_{r}(H_f).
$$
 By the definition of $\Pi_{p, 2}(F^p_\alpha, L^p_\alpha)$,  for any fixed integer $m$, we see
\[
\left( \sum_{j=1}^{m} \| H_f(k_{z_j})\|_{p, \alpha}^p \right)^{\frac{1}{p}} \lesssim \pi_{p,2}(H_f) \cdot \sup_{g \in B_{F_\alpha^{p'}}} \left( \sum_{j=1}^{m} |\langle g, k_{z_j} \rangle|^2 \right)^{\frac{1}{2}},
\]
where $\{z_j\}$ is a $\delta$-lattice with $\delta > 0$ small enough. Since $0< p'/2\leq 1$ for $p\geq 2$,  we know
\begin{align*}
\left( \sum\limits_{j=1}^{m}\|H_f(k_{z_j})\|_{p, \alpha}^{p}\right)^{\frac{p'}{p}}  &\leq \pi_{p,2}(H_f)^{p'} \cdot\left[ \sup_{h\in B_{F^{p'}_\alpha}} \left(\sum\limits_{j=1}^{m}\left| \langle h, k_{z_j} \rangle \right|^{2}\right)^{\frac{1}{2}}\right]^{p'}  \\
&\leq \pi_{r}(H_f)^{p'} \cdot \sup_{h\in B_{F^{p'}_\alpha}} \sum\limits_{j=1}^{m}\left| \langle h, k_{z_j} \rangle \right|^{p'} \\
&\simeq\pi_{r}(H_f)^{p'} \cdot \sup_{h\in B_{F^{p'}_\alpha}} \sum\limits_{j=1}^{m}\left|h(z_j)e^{-\frac{\alpha}{2}|z_j|^2}\right|^{p'}\\
&\lesssim \pi_{r}(H_f)^{p'}.
\end{align*}
On the other hand, similarly to the proof in Case 1, we have
\begin{align*}
G_{p,2\delta}(f)(z_j)^{p} \lesssim  \|H_f(k_{z_j})\|_{p, \alpha}^{p}.
\end{align*}
Hence,
  \begin{align*}
\sum\limits_{j=1}^{m}G_{p, 2\delta}(f)(z_j)^{p} \lesssim \sum\limits_{j=1}^{m}\|H_f(k_{z_j})\|_{p, \alpha}^{p}
\lesssim \pi_{r}(H_f)^{p} .
\end{align*}
Letting $m\rightarrow\infty$, we obtain
 \begin{align*}
\int_{{\mathbb C}^n}G_p(f)^{p}dv\lesssim \sum\limits_{j=1}^{\infty}G_{p, 2\delta}(f)(z_j)^{p}\lesssim \pi_{r}(H_f)^p<\infty.
\end{align*}
\end{proof}

\section{ Berger-Coburn phenomenon}\label{sect4}

In this section, we examine the Berger-Coburn phenomenon (BCP) for absolutely summing Hankel operators on Fock spaces. BCP refers to the property where, for a bounded symbol \( f \), the Hankel operator \( H_f \) belongs to a given operator ideal \( \mathcal{M} \) if and only if  \( H_{\bar{f}} \) also belongs to \( \mathcal{M} \). Formally, for an operator ideal \( \mathcal{M} \), BCP holds if
\[
H_f \in \mathcal{M} \iff H_{\bar{f}} \in \mathcal{M}, \quad \forall f \in L^\infty,
\]
see \cite{HV22} for details.

BCP were identified by Berger and Coburn in \cite{BC87}, Hagger and   Virtanen in \cite{HaV21} and Zhu in \cite{Zh12} for compact operators on the classical Fock space.  On schatten classes, BCP were  discussed by Bauer in \cite{Ba04},  Xia and Zheng in \cite{XZ}, Hu and Virtanen in \cite{HV22},   Asghari, Hu, and Virtanen  in \cite{AHV24}.   Xia \cite{Xia23} provided a counterexample of some function in $L^\infty$ showing that \( H_f \) can be trace-class while \( H_{\bar{f}} \) is not.

In the context of BCP for absolutely summing   operators  ideals, we need Proposition \ref{por01}.  It presents a key norm estimate for functions in the space  $\mathrm{IDA}^{s, p}$, which we believe that it has its own interest.

\begin{proposition}\label{por01}
Suppose \( 1 < s < \infty , 1 \leq p < \infty \).  For  \(f \in L^{\infty}(\mathbb{C}^{n})\cap \mathrm{IDA}^{s,p} \), there holds
\[
\|\overline{f}\|_{\mathrm{IDA}^{s,p}}\lesssim\|f\|_{\mathrm{IDA}^{s,p}}.
\]
\end{proposition}

\begin{proof} For \( f \in L^{\infty}(\mathbb{C}^{n}) \) with \(\|f\|_{\mathrm{IDA}^{s,p}} < \infty\),  \cite[Lemma 3.6]{HV23} gives the decomposition \( f = f_1 + f_2 \) as
\[
f_1 = \sum_{j=1}^{\infty} h_j \psi_j\in C^2(\mathbb{C}^{n}) \quad \textrm{and} \quad f_2 = f - f_1,
\]
where $h_j$ and $\psi_j$ defined as in Lemma 3.5 in \cite{HV23}. Moreover, for any $z\in \mathbb{C}^{n}$, we have
 \[
\left| \overline{\partial}f_1(z) \right| + M_{p, r/2} (\overline{\partial}f_1)(z) + M_{p, r/2} (f_2)(z)
\lesssim G_{p, 2r}(f)(z).
\]
 For $1\leq p<\infty$,   by the estimate (3-20) of \cite{HV23} and \cite[Corollary 3.10]{HV23}, we can choose some \(\tau> 0 \) such that
\[
\|M_{p, r}(\overline{\partial}f_1)\|_{L^s} + \|M_{p, r}(f_2)\|_{L^s} \lesssim \|G_{p, \tau}(f)\|_{L^s}\simeq  \|f\|_{\mathrm{IDA}^{s,p}}.
\]
Hence,
\[
\|M_{p,r}(\overline{\partial}f_1)\|_{L^s} \lesssim \|f\|_{\mathrm{IDA}^{s,p}}.
\]
Notice that, \( f_1 \in C^2 \cap L^{\infty} \). Similarly to the proof in \cite{HV22}, we   claim
\begin{equation}\label{por02}
\left(\int_{\mathbb{C}^{n}}M_{p, r/4}(\overline{\partial}  \, \overline{f_1})^sdv\right)^{1/s}\lesssim  \|f\|_{\mathrm{IDA}^{s,p}}.
\end{equation}
In fact, as on \cite[Page 6012]{HV23a}, there exist $F$ and $H$ such that \( \overline{\partial} \, \overline{f_1} = F + H \), where
\[
F = \sum_{j=1}^{\infty} \overline{h_j} \cdot \overline{\partial} \, \psi_j , \quad H = \sum_{j=1}^{\infty} \psi_j \cdot\overline{\partial} \, \overline{h_j}.
\]
With the same proof as the one for the estimate of $\left|\overline{\partial}f_1\right|$ in \cite{HV22}, we have
\[
|F(z)| \lesssim G_{p, 2r}(f)(z).
\]
For $1 < s < \infty$, by \cite[Lemma 7.1]{HV22}, we obtain
\[
\|H\|_{L^s} \leq \|\overline{\partial} \, \overline{f_1}\|_{L^s} + \|F\|_{L^s} \lesssim  \|\overline{\partial}f_1\|_{L^s} + \|F\|_{L^s}
 \lesssim  \|f\|_{\mathrm{IDA}^{s,p}}.
\]
As the proof in \cite{HV23a}, we obtain \eqref{por02}.
By \cite[Theorem 3.8]{HV23}, we have
 \begin{equation}\label{anorm01}
M_{p, r/4}(\overline{\partial}f_1)+ M_{p, r/4}(f_2) \in L^s.
 \end{equation}
Furthermore, \[
\|f\|_{\mathrm{IDA}^{s, p}}\simeq \inf\left\{\left\|M_{p, \, r/4}(\overline{\partial}f_1)\right\|_{L^s}+\left\|M_{p, \, r/4}(f_2)\right\|_{L^s}\right\},
\]
where the infimum is taken over all possible decompositions $f = f_1 + f_2$ that satisfy \eqref{anorm01}. Since
$$
\left\|M_{p, \, r/4}(f_2)\right\|_{L^s}=\left\|M_{p, \, r/4}(\overline{f_2})\right\|_{L^s},
$$
\eqref{por02} and  \cite[Corollary 3.10]{HV23} yield
\[
\| \overline{f} \|_{\mathrm{IDA}^{s,p}} \leq \left\|M_{p, \, r/4}(\overline{\partial} \, \overline{f_1})\right\|_{L^s}+\left\|M_{p, \, r/4}(\overline{f_2})\right\|_{L^s}\lesssim \|f\|_{\mathrm{IDA}^{s,p}}.
\]
\end{proof}

{\bf Remark}:  If \(0 < s \leq 1\), the conclusion of  Proposition \ref{por01} need not be true. To see this we  consider the case where \( n = 1 \), and define $f(z)$ as
\[
f(z) =
\begin
{cases}
\frac 1 z, &
\text{if } |z| \geq 1; \\
0, &
\text
{if } |z| < 1.
\end
{cases}
\]
An  elementary   computation shows $f\in \mathrm{IDA}^{s,p}$, while $\overline{f}\notin \mathrm{IDA}^{s,p}$ for all $p\geq 1$.\\

Notice that, $\kappa>1$ as in \eqref{number-s1} if $1\leq p<\infty$ and  $r\geq1$. Combining Theorem \ref{thm:main} and Proposition \ref{por01} above, we have the following theorem.

\begin{theorem}\label{thm:main-5}
Suppose $1\leq p<\infty$, $r\geq1$. Then  the BCP holds for $\Pi_r(F^p_\alpha, L^p_\alpha)$, that is, for $f \in L^\infty(\mathbb{C}^{n})$, $H_f\in \Pi_r(F^p_\alpha, L^p_\alpha)$
if and only if $H_{\overline{f}}\in \Pi_r(F^p_\alpha, L^p_\alpha)$ with the norm estimate
$$
\pi_r(H_f)\simeq \pi_r(H_{\overline{f}}).
$$
\end{theorem}

\section{final remark}

It is clear from the proof that our results also hold for the weighted Fock spaces $F^p_\varphi$ with   $\varphi$ satisfying condition \eqref{a0bwei}. Meanwhile, as introduced in \cite{Ch91, MMO03}, the doubling  Fock space is certain generalization of $F^p_\alpha$ in the case $n=1$. We believe that, with appropriate modifications to the proof, our method can be extended to the doubling Fock spaces.

\end{document}